%

\input amstex
\define\qued{~\qed \enddemo}

\define\forces{\Vdash}
\define\nforces{\nVdash}
\define\clinof#1#2{\text{cl}_{#1}{#2}}
\define\eset{E^\omega_{\omega_2}}
\define\cl#1{\overline{#1}}
\define\etheta{E^\omega_\theta}
\define\fsquare{F_{\theta,\omega}\times F_{\theta,\omega}}

\documentstyle{amsppt}
\loadbold
\magnification=1200
\NoRunningHeads

\topmatter
\title \nofrills
{{\bf REFLECTION AND WEAKLY \\ COLLECTIONWISE HAUSDORFF SPACES}} \\
 \\  {\rm Tim LaBerge}\\{\rm Avner Landver}
\endtitle
\address
Department of Mathematics,
The University of Kansas,
Lawrence, KS 66045
\endaddress

\date October 20, 1992 \enddate

\email
LaBerge\@kuhub.cc.ukans.edu Landver\@kuhub.cc.ukans.edu
\endemail

\abstract 
We show that $\square(\theta)$ implies that there is a first
countable $<\!\theta$-collectionwise Hausdorff space that is not weakly
$\theta$-collectionwise Hausdorff.  We also show that in the model obtained 
by Levy collapsing a weakly compact (supercompact) cardinal to $\omega_2$,
first countable $\aleph_1$-collectionwise Hausdorff spaces are weakly
$\aleph_2$-collectionwise Hausdorff (weakly collectionwise
Hausdorff).  In the last section we show that assuming $E^\omega_\theta$,
a certain $\theta$-family of integer valued functions exists and that
in the model obtained by Levy collapsing a supercompact
cardinal to $\omega_2$, these families do not exist.  
\endabstract



\endtopmatter

\document

\subhead \S 1 Introduction \endsubhead
Reflection is a central theme in modern set-theoretic topology. As Alan Dow
points out in \cite{Do}, we often prove theorems when some type of
reflection principle holds, and we build counter-examples when
reflection fails.  
This paper contains both types of results, on questions related to the
failure of collectionwise Hausdorff. 

We say that a subset $A$ of a
topological space $X$ can be {\it separated} if there is a collection $\{U_x
: x \in A\}$ of disjoint open sets with $x \in U_x$ for every $x \in A$.   A
space $X$ is {\it $<\!\theta$-collectionwise Hausdorff} ($<\!\theta$-cwH) if
every closed discrete subset of size  $<\!\theta$ can be separated.  A
space $X$ is {\it $\le\!\theta$-collectionwise Hausdorff}
($\le\!\theta$-cwH) if every closed discrete subset of size  $\le\theta$ can
be separated. $X$ is
{\it collectionwise Hausdorff} (cwH) if it is $\le\!\theta$-cwH for every
cardinal $\theta$. 

In particular, we are concerned with Fleissner's questions:

\proclaim {Question 1} Is there a ZFC example of a first countable 
$<\!\aleph_2$-cwH space that is not $\le\!\aleph_2$-cwH?
\endproclaim 

\proclaim {Question 2} Is there a ZFC example of a first countable 
$<\!\aleph_2$-cwH space that is not cwH?
\endproclaim

Fleissner asks for ZFC examples because he  showed \cite{F} that $\eset$
(a non-reflecting stationary subset of $\{\alpha<\omega_2: cf(\alpha) =
\omega\}$) can be used to construct an example of a first countable 
$<\!\aleph_2$-cwH space that is not $\le\!\aleph_2$-cwH.  Since the failure
of $\eset$ is equiconsistent with the existence of a Mahlo cardinal
(\cite{De} and \cite{HS}),  large cardinals are required for a negative
answer to Question~1.   In fact, we'll see that by a result of Todorcevic
\cite{To2}, at least a weakly compact cardinal will be needed to get a
negative answer to Question~1.  Shelah has shown \cite{S} that in the model
obtained by Levy collapsing a weakly compact cardinal to $\omega_2$, first
countable $<\!\aleph_2$-cwH spaces that are locally of size at most
$\aleph_1$ are $\le\!\aleph_2$-cwH, and that in the model obtained by
collapsing a supercompact cardinal, such spaces are cwH. For a more
complete history of the general problem of reflecting the failure of cwH,
see \cite{W} and \cite{FS}.  We would like to thank Bill Fleissner for his
valuable input and  Gary Gruenhage for letting us use the
handwritten notes \cite{GT}.

\subhead 
\S 2 Squares of Fans and First Countable Collectionwise Hausdorff Spaces 
\endsubhead

Recently, Gruenhage and Tamano \cite{GT} have discovered a connection
between Fleissner's questions and the problem of determining the tightness of
the squares of certain fans.  Before we discuss this connection, we make some
definitions.

The fan $F_{\theta,\omega}$ is the quotient space obtained by identifying the
non-isolated points of $\theta$-many copies of the convergent
$\omega$-sequence.  To be precise, $F_{\theta,\omega} =
\{*\}\cup(\theta\times\omega)$, topologized so that points of
$\theta\times\omega$ are isolated and so that an  open  base at
$*$ is the family of all 
$$
B_f = \{*\}\cup\{(\alpha,n):n\ge f(\alpha)\} \quad (f\in{}^\theta\omega).
$$

When working with the square of $F_{\theta,\omega}$, one can always replace
two functions $f_1,f_2\in{}^\theta\omega$ with $g = \max(f_1,f_2)$, so we use
the family of all 
$$
V_g = B_g\times B_g\qquad (g\in{}^\theta\omega)
$$
as an open base for  $(*,*)$ in $\fsquare$.

We say that a set $S\subseteq \fsquare$ is
{\it $\theta$-good} if $(*,*)\in\cl{S}$, but for all $T\in [S]^{<\theta}$,
$(*,*)\notin \cl{T}$.  Thus, ``there is a $\theta$-good set in
$\fsquare$'' means that the tightness of
$\fsquare$ is $\theta$ and that this
tightness is actually attained.

Finally, we say that a set $H\in[\omega\times\omega]^{<\omega}$ is {\it
closed downward} (c.d.w.) if whenever $(n,m)\in H$, then $n'\le
n$ and $m'\le m$ implies $(n',m')\in H$.

The Gruenhage-Tamano result is the equivalence (1)$\iff$(2) in
Theorem~1.  We have added the combinatorial equivalence~(3) and the 
superficially stronger topological characterization~(4).

\proclaim {Theorem~1}
Let $\theta$ be an uncountable cardinal.  TFAE
\roster
\item  There exists a $\theta$-good subset of $\fsquare$.

\item  There exists a space $X$ that is first countable, $<\!\theta$-cwH and
not $\le\!\theta$-cwH.

\item  There exists ${\Cal H}=\{H_{\alpha\beta} : \alpha < \beta < \theta\} 
\subseteq [\omega \times \omega]^{<\omega}$ with each 
$H_{\alpha\beta}$ c.d.w. such that

\itemitem{(a)}  For every $A \in [\theta]^{<\theta}$ there is a function
$f : \theta \to \omega$ such that for all $\alpha < \beta$  in~$A$,
$(f(\alpha),f(\beta)) \notin H_{\alpha\beta}$.

\itemitem{(b)} For every $f:\theta \to \omega$ there are $\alpha <\beta <
\theta$ such that $(f(\alpha),f(\beta))\in H_{\alpha\beta}$.

\item  There exists a space $X$ that is first countable, zero-dimensional,
$<\!\theta$-cwH and not $\le\!\theta$-cwH.
\endroster
\endproclaim

In order to put Theorem~1 in perspective, we need the definition of a
difficult-to-deny combinatorial principle.  For $\theta$   an
uncountable regular cardinal, $\square(\theta)$ is the assertion that there is
a family  $\{C_\alpha : \alpha < \theta\}$ satisfying the following conditions:
\roster
\item "(i)" $C_\alpha \subset \alpha$ is a club subset of $\alpha$.
\item "(ii)" If $\alpha$ is a limit point of $C_\beta$, then 
$C_\alpha = C_\beta \cap \alpha$.
\item "(iii)" There is no club $C \subset \theta$ such that for every
   limit point $\alpha$ of $C$, $C_\alpha = C \cap \alpha$.
\endroster
$\square(\theta)$ is true for every regular $\theta$ which is not weakly
compact in $L$ (see \cite{To1}).

Todorcevic showed (see \cite{To2} and \cite{Be}) that $\square(\theta)$
can  be used to construct a $\theta$-good set.  Also, combining Theorem~1
with Fleissner's construction of a first countable, $<\!\theta$-cwH, not
$\le\!\theta$-cwH space from $\etheta$, we see that $\etheta$ 
can be used to construct a $\theta$-good set.

Because the proof of Theorem~1 is essentially the same as the proof of
Theorem~2, we only prove Theorem~2.
Before the statement of the theorem we  define a weakening of 
collectionwise Hausdorff introduced by Tall \cite{Ta}. We say that a subset
$A$ of a space $X$ is {\it weakly separated\/} if it has a subset of size
$|A|$ that is separated.   $X$ is
{\it weakly $\theta$-cwH\/} if every closed discrete subset of size
$\theta$ is weakly  separated.  $X$ is {\it weakly cwH\/} if it is weakly
$\theta$-cwH for every cardinal $\theta$.

\proclaim {Theorem~2}
Let $\kappa \le \theta$ be uncountable cardinals.  TFAE
\roster

\item  There exists a space $X$ that is first countable, $<\!\kappa$-cwH and
not weakly $\theta$-cwH.

\item  There exists ${\Cal H}=\{H_{\alpha\beta} : \alpha < \beta < \theta\} 
\subset [\omega \times \omega]^{<\omega}$, with each $H_{\alpha\beta}$  c.d.w.,
such that

\itemitem{(a)}  For every $A \in [\theta]^{<\kappa}$ there is a function
$f : \theta \to \omega$
such that for all $\alpha < \beta$  in~$A$, $(f(\alpha),f(\beta)) \notin
H_{\alpha\beta}$.

\itemitem{(b)}  For every $B \in [\theta]^{\theta}$ and every
$f:\theta \to \omega$ there are $\alpha <\beta$ in $B$
such that  $(f(\alpha),f(\beta))\in H_{\alpha\beta}$.

\item  There exists a space $X$ that is first countable, zero-dimensional,
$<\!\kappa$-cwH and not weakly $\theta$-cwH.
\endroster
\endproclaim

\demo {Proof}
$(1) \Rightarrow (2)$:  Let $\theta \subset X$ be a subset which is not
weakly separated.  For every $\alpha \in \theta$ let $\{U_n(\alpha) : n \in
\omega \}$ be a decreasing neighbourhood base at $\alpha$.  Let 
$$
S =\{ (\langle\alpha, n\rangle,\langle\beta, m\rangle) : 
U_n(\alpha) \cap U_m(\beta) \ne \emptyset \}.
$$
$S$ is a subset of $\fsquare$.   
For every $B \subset \theta$, let  
$S_B =\{ (\langle\alpha, n\rangle,\langle\beta, m\rangle) \in S : \alpha,
\beta \in B \}$. We claim that $S$ satisfies the following two conditions 
\roster
\item "(i)"  For every $A \in [\theta]^{<\kappa}$, $(*,*) \notin \cl S_A$.
\item "(ii)"  For every $B \in [\theta]^{\theta}$, $(*,*) \in \cl S_B$.
\endroster

To show (i), let $A \in [\theta]^{<\kappa}$.  Let $f:A \to \omega$ be such
that $\{U_{f(\alpha)}(\alpha) : \alpha \in A \}$ is a separation of $A$ in
$X$. Let $V_f$ be the open neighbourhood of $(*,*)$ in
$\fsquare$ that is determined by $f$. It is
easy to check that $V_f \cap S_A = \emptyset$. For (ii), let $B \in
[\theta]^{\theta}$. Let $f : \theta \to \omega$, and let $V_f$ be as
above.   Now, $f \!\restriction\! B$ is not a code for a separation of $B$ in
$X$, therefore there are $\alpha < \beta$ in $B$ such that 
$U_{f(\alpha)}(\alpha) \cap U_{f(\beta)}(\beta) \ne \emptyset$.  But then
$(\langle\alpha,f(\alpha)\rangle, \langle\beta,f(\beta)\rangle) \in S_B \cap
V_f$.

Using the fact that $(*,*)$ is not in the closure of any countable subset
of $S$, it is not hard to verify the following facts.

\proclaim {Fact 1}
For every $\alpha < \theta$ there is $h(\alpha) \in \omega$ such that for
every $n > h(\alpha)$ and for every $\beta < \theta$  the set $\{m\in \omega:
(\langle \alpha, n\rangle,\langle\beta, m\rangle) \in S \}$ is finite.
\endproclaim

\proclaim {Fact 2}
For every $\beta < \theta$ there is $g(\beta) \in \omega$ such that for
every $m > g(\beta)$ and for every $\alpha < \theta$  the set $\{n\in
\omega:  (\langle \alpha, n\rangle,\langle\beta, m\rangle) \in S \}$ is
finite. \endproclaim

Let $T = S \cap (B_h \times B_g)$.  Clearly, $T$ satisfies~(i) and~(ii).
Now, for every $\alpha < \beta < \theta$ define
$$
H'_{\alpha\beta}=\{(n,m)\in \omega\times\omega :  
(\langle \alpha, n\rangle,\langle\beta, m\rangle) \in T\}.
$$

Let $H_{\alpha\beta}$ be the downward closure of $H'_{\alpha\beta}$. Let us show
that the $H_{\alpha\beta}$'s are finite.  It is enough to show that the
$H'_{\alpha\beta}$'s are finite.  Assume that this is false and let $\alpha
< \beta < \theta$ be such that $H'_{\alpha\beta}$ is infinite. It follows
that for every $k \in \omega$ there are $n, m > k$ such that $(n, m) \in
H'_{\alpha\beta}$.  This implies that  $(*, *)$   is in the closure of
$\{(\langle\alpha, n\rangle, \langle\beta ,m\rangle): (n, m) \in H'_{\alpha\beta}
\}$, therefore $(*, *)$ is in the closure of a countable subset of $T$,  a
contradiction.  

Finally, it is not hard to see that 
$\{H_{\alpha\beta} : \alpha < \beta < \theta \}$ satisfies (a) and (b) of (2)
precisely because $T$ satisfies (i) and (ii).

$(2) \Rightarrow (3)$:  As in \cite{GT}, first let  $I =\{ (\langle \alpha,
n\rangle,\langle\beta, m\rangle) : \alpha < \beta < \theta \land (n,m)  \in
H_{\alpha\beta} \}$.  Then let $X=I \cup \theta$.  Points in $I$ are 
isolated, and for every $\gamma \in \theta$ a  base at $\gamma$ is 
given by 
$$
U_k (\gamma) = \{\gamma\} \cup \{ (\langle \alpha, n\rangle,\langle\beta,
m\rangle) \in I:  (\alpha = \gamma \land n \ge k) \lor (\beta = \gamma \land
m \ge k) \} \qquad(k\in\omega).
$$

Clearly, $X$ is first countable.  To see that $X$ is $<\!\kappa$-cwH let
$A \in [\theta]^{<\kappa}$.  Let $f : \theta \to \omega$ be given by (2a).
Then check that $\{U_{f(\alpha)}(\alpha): \alpha \in A \}$ is a separation
of $A$.  To see that $X$ is not weakly $\theta$-cwH let  $B \in
[\theta]^{\theta}$ and let  $f : \theta \to \omega$.  By (2b), there are 
$\alpha <\beta$ in $B$ such that $(f(\alpha),f(\beta))\in H_{\alpha\beta}$.  So,
$$
(\langle\alpha,f(\alpha)\rangle,\langle\beta, f(\beta)\rangle) \in
U_{f(\alpha)}(\alpha) \cap U_{f(\beta)}(\beta),
$$
therefore $f$ does not code a separation of $B$.

Finally, the finiteness of the $H_{\alpha\beta}$'s implies that
each $U_k(\alpha)$ is clopen and therefore $X$ is zero-dimensional.

$(3) \Rightarrow (1)$: Trivial.\qued

\subhead \S 3 Independence Results \endsubhead

In this section, we first show that $\square(\theta)$ implies that there is a
first countable, $<\!\theta$-cwH space that is not weakly $\theta$-cwH.  We
then demonstrate the consistency of ``first countable $<\!\aleph_2$-cwH spaces
are weakly $\aleph_2$-cwH (weakly cwH)'', assuming the consistency of ``there
is a weakly compact (supercompact) cardinal''.

Since $\square(\theta)$ is true unless $\theta$ is weakly compact in $L$, this
gives that ``first countable $<\!\aleph_2$-cwH spaces
are weakly $\aleph_2$-cwH'' is equiconsistent with ``there is a weakly compact
cardinal'' and that the consistency of ``first countable $<\!\aleph_2$-cwH spaces
are weakly cwH'' implies the existence of an inner model with many measurable
cardinals~\cite{KM}.

\proclaim{Theorem~3} Let $\theta$ be a regular uncountable cardinal and suppose
that $\square(\theta)$ holds.  Then there is a zero-dimensional, first
countable, $<\!\theta$-cwH  space that is not weakly $\theta$-cwH. \endproclaim

\demo{Proof}  In \cite{Be}, Todorcevic constructs, from the assumption of
$\square(\theta)$, a   function  $\rho_2$ that maps pairs
$\alpha<\beta< \theta$ into $\omega$. In $\S 4$, we will study
integer-valued families of functions, so to keep our notation consistent,
we define for each $\beta<\theta$ a function $h_\beta:\beta\to \omega$ by
$$
h_\beta(\alpha) = \rho_2(\alpha,\beta)\qquad(\alpha<\beta).
$$
By the properties of $\rho_2$ cited in \cite{Be}, the family
$\{h_\beta:\beta<\theta\}$ satisfies
\roster 
\item For all $\alpha<\beta< \theta$, there is an
$n_{\alpha\beta}\in\omega$ such that for all $\xi<\alpha$,
$h_\alpha(\xi)<h_\beta(\xi)+n_{\alpha\beta}.$ (This is monotonicity
\cite{DW}.)
\item For all $B\in[\theta]^\theta$ and all $n\in\omega$, there
exist $\alpha<\beta$ in $B$ such that $h_\beta(\alpha)>n.$ \endroster

For every $\alpha<\beta< \theta$, define 
$H_{\alpha\beta} = \{(n,m): n+m\le h_\beta(\alpha)\}$.
Let us show that $\Cal H =\{H_{\alpha\beta}:\alpha<\beta< \theta\}$ 
satisfies~(2) of Theorem~2 with $\kappa =\theta$.

For part (a) let $A \in [\theta]^{<\kappa}$.  Let $\gamma < \theta$ be a bound
for $A$ and define  $f : \theta \to \omega$ by
$$
f(\alpha) = \cases 
h_\gamma(\alpha) + n_{\alpha\gamma}
&\text{if } \alpha < \gamma, \\
0 & \text {otherwise}.
\endcases
$$
Now let $\alpha < \beta$ in $A$, then $f(\alpha)+f(\beta) \ge 
h_\gamma(\alpha) + n_{\beta\gamma} > h_\beta(\alpha)$.  Therefore 
$(f(\alpha),f(\beta)) \notin H_{\alpha\beta}$.

For part (b) let $B \in [\theta]^{\theta}$ and $f:\theta \to \omega$.
Then there exist $B' \in [B]^\theta$ and $n \in \omega$ such that for every
$\alpha \in B'$, $f(\alpha)=n$.  But then there are $\alpha < \beta$ in $B'$
such that $h_\beta(\alpha) > 2n$.  Hence 
$f(\alpha) + f(\beta) = 2n < h_\beta(\alpha)$ and so 
$(f(\alpha),f(\beta)) \in H_{\alpha\beta}$.
\qued

\remark{Remark}  Since the $h_\beta$'s also satisfy

\noindent (3) Whenever $A,B\in[\theta]^\theta$ and $n\in\omega$, there are
$\alpha<\beta$ with $\alpha\in A$ and $\beta\in B$ such that
$h_\beta(\alpha)>n$, 

\noindent it follows that the space constructed above is badly non-normal. 
\endremark

The next two theorems will use the method of forcing $+$ reflection.  A good
reference for this technique is \cite{DTW1} and \cite{DTW2}.  We first prove
a preservation lemma, whose proof is  motivated by Lemma 3.13 in \cite{Be}.

\proclaim {Main Lemma}
Assume that ${\Bbb P}$ is an $\omega_1$-closed partial order 
and that $\theta \ge \omega_2$ is a regular cardinal.  Assume that 
${\Cal H}=\{H_{\alpha\beta} : \alpha < \beta < \theta\} 
\subset [\omega \times \omega]^{<\omega}$ are c.d.w.\ and the following 
condition from Theorem~2 holds:
\roster
\item "(b)" For every $B \in [\theta]^{\theta}$ and every
$f:\theta \to \omega$ there are $\alpha <\beta$ in $B$
such that $(f(\alpha),f(\beta))\in H_{\alpha\beta}$.
\endroster
Assume that $G$ is a ${\Bbb P}$-generic filter over $V$.  Then in
$V[G]$, for every $f : \theta \to \omega$ there are $\alpha < \beta <
\theta$ such that  $(f(\alpha),f(\beta)) \in H_{\alpha\beta}$. 
\endproclaim

\demo {Proof}
Assume that the lemma is false.  Then there is $p_0 \in {\Bbb P}$ and 
 a ${\Bbb P}$-name $\dot g$ such that 
$$
p_0 \forces  ``\dot g: \theta \to \omega \land  (\forall \alpha < \beta
< \theta)\, (\dot g(\alpha), \dot g(\beta)) \notin H_{\alpha \beta}".
$$

Let $M$ be an elementary substructure of some $H(\lambda)$ ($\lambda$ large
enough) with  $\theta$, ${\Cal H}$, ${\Bbb P}$, $p_0$, $\dot g$, and 
$\forces_{\Bbb P}$ in $M$ and   $M \cap \theta = \beta < \theta$.

If possible, choose $p_1 \le p_0$ and $\alpha_0 < \beta$ such that $p_1 \in
M$ and $p_1 \Vdash ``(\dot g(\alpha_0), 0) \in H_{{\alpha_0} \beta}"$. Note
that this would imply that $p_1 \Vdash ``\dot g(\beta) > 0"$.

Similarly construct $p_0 \ge p_1 \ge \dots \ge p_{n+1} \dots$ such that
$p_{n+1} \in M$ and $\alpha_n < \beta$ with
$$
p_{n+1} \forces `` (\dot g(\alpha_n), n) \in H_{{\alpha_n} \beta}".
$$

In particular this would imply that $p_{n+1} \Vdash `` \dot g(\beta) > n"$.
This process must stop at a finite stage since otherwise one could choose 
$q$ with $q \le p_n$ for every $n \in \omega$ and get that
for every $n \in \omega$, $ q \forces `` \dot g(\beta) > n"$, 
which is impossible.
\qued

So, there is $\bar n \in \omega$ such that for every $q \le p_{\bar n}$ and
for every $\alpha < \beta$, if $q \in M$, then 
$$
q \nforces ``(\dot g(\alpha), \bar n) \in H_{\alpha \beta}".
$$

\proclaim {Claim 1}
For every $\alpha < \beta$, 
$p_{\bar n} \forces ``(\dot g(\alpha), \bar n) \notin H_{\alpha \beta}"$.  
\endproclaim

\demo {Proof}
Assume not.  Then there is $q \le p_{\bar n}$ and $\alpha < \beta$  such that 
$ q \forces ``(\dot g(\alpha), \bar n) \in H_{\alpha \beta}"$.  Fix such 
$\alpha < \beta$ and the corresponding $H_{\alpha \beta}$.
Then $H(\lambda) \models (\exists q \le p_{\bar n})\,\, q \forces 
``(\dot g(\alpha), \bar n) \in H_{\alpha \beta}"$.
But $p_{\bar n}$, $\Vdash_{\Bbb P}$, $\dot g$, $\alpha$ and $H_{\alpha \beta}$
are all in $M$, therefore  
$$
M \models    \exists q \le p_{\bar n} 
\,\,  (q \forces ``(\dot g(\alpha), \bar n)
\in H_{\alpha \beta}"),
$$
which is a contradiction.
\qued

Let $B=\{ \beta \in \theta :  \forall \alpha < \beta  \,\,
(p_{\bar n} \forces ``(\dot g(\alpha), \bar n) \notin H_{\alpha \beta}")\}$. 
By elementarity of $M$, $|B|=\theta$.  

Now, for every $\alpha \in \theta$ let $q_\alpha \le p_{\bar n}$ and 
$m_\alpha \in \omega$ be such that $q_\alpha \forces ``\dot{g}(\alpha) =
m_\alpha"$. In $V$, define a function $f : \theta \to \omega$ by:
$$
f(\alpha) = m_\alpha + {\bar n}\qquad (\alpha\in\theta).
$$

The following claim will give us the desired contradiction.

\proclaim {Claim 2}
For every $\beta \in B$ and every $\alpha < \beta$, $(f(\alpha),f(\beta)) \notin
H_{\alpha \beta}$. \endproclaim

\demo {Proof}
Fix $\beta \in B$ and  $\alpha < \beta$.  By definition,
$f(\alpha) \ge m_\alpha$ and $f(\beta) \ge {\bar n}$, therefore it is enough 
to show that $ (m_\alpha, {\bar n}) \notin H_{\alpha \beta}$.  
But $p_{\bar n} \forces ``(\dot g(\alpha), \bar n) \notin H_{\alpha \beta}"$,
and $q_\alpha \le p_{\bar n}$, and 
$q_\alpha \forces `` \dot g(\alpha) = m_\alpha"$.
Therefore $q_\alpha \forces ``(m_\alpha,\bar n) \notin H_{\alpha \beta}"$,
and this implies that indeed $(m_\alpha,\bar n) \notin H_{\alpha \beta}$.
\qued

Let $\kappa$ be strongly
inaccessible and let ${\Bbb P}={\Bbb P_\kappa}$ be the Levy collapse of
$\kappa$ to  $\omega_2$ with countable conditions 
(for a definition and proofs  of the facts below, see \cite{K}). Let $G$
be a  ${\Bbb P}$-generic filter. We will use the following facts about
${\Bbb P}$:  ${\Bbb P}$ is an  $\omega_1$-closed partial order.  ${\Bbb P}$
is $\kappa$-c\.c\.  and $\kappa=\omega_2$ in $V[G]$. For every $\lambda <
\kappa$, ${\Bbb P}$ can be factored as  ${\Bbb P} = {\Bbb P}_\lambda \times
{\Bbb P}^\lambda$, ${\Bbb P} \cap V_\lambda = {\Bbb P}_\lambda$ and in
$V[G_\lambda]$,  ${\Bbb P}^\lambda$ is forcing equivalent to ${\Bbb P}$ (in
particular it is $\omega_1$-closed). For a definition and discussion of
weakly compact cardinals see \cite{K} and \cite{KM}.  The following type
of proof was invented by Baumgartner \cite{Ba}, see also \cite{DTW2}.

\proclaim{Theorem~4} Assume that $\kappa$ is a weakly compact cardinal.  
Let ${\Bbb P}$ be the Levy collapse of $\kappa$ to $\omega_2$ with countable
conditions and let $G$ be a 
${\Bbb P}$-generic filter over $V$.  Then in $V[G]$, every first
countable $<\!\aleph_2$-cwH space is weakly $\aleph_2$-cwH. 
\endproclaim

\demo{Proof}  
By way of contradiction, suppose that in $V[G]$ there is a  first
countable $<\!\aleph_2$-cwH space that is not weakly $\aleph_2$-cwH.  In
$V[G]$, $\omega_2 = \kappa$, therefore there is  ${\Cal H} =
\{H_{\alpha\beta}:\alpha<\beta< \kappa\}$ that satisfies  $(2)$ of
Theorem~2 with $\kappa=\theta$. This fact must be forced and the forcing
statement is $\Pi_1^1$ over $(V_\kappa, \in)$ (with few extra parameters).
Using the facts that $\kappa$ is weakly compact ($\Pi_1^1$-reflection) 
and that ${\Bbb P}$ is $\kappa$-c\.c\., we find an inaccessible 
$\lambda<\kappa$ such that in $V[G_\lambda]$, 
$\{H_{\alpha\beta}:\alpha<\beta<\lambda\}$ satisfies $(2)$ of Theorem~2
(with $\lambda$ playing the role of both $\kappa$ and $\theta$).

Now, $V[G_\lambda] \models  {{\Bbb P}^\lambda}$ is
$\omega_1$-closed.  Therefore, by the Main Lemma, for  
every function ${f}:\lambda\to \omega$ in $V[G]$, there are 
$\alpha<\beta<\lambda$ such that $(f(\alpha), f(\beta))\in H_{\alpha\beta}$. 
But $\lambda < \kappa$, therefore
this contradicts the fact that (a) of part~(2) of Theorem~2 holds for 
$\Cal H$ in $V[G]$.\qued

For a definition and discussion of supercompact cardinals see
\cite{KM}, \cite{J}, and \cite{DTW1}.

\proclaim {Theorem~5}
Assume that $\kappa$ is a supercompact cardinal and that the GCH holds above
$\kappa$.  Let ${\Bbb P}$ be the Levy collapse of $\kappa$ to $\omega_2$ with
countable conditions and let $G$ be a ${\Bbb P}$-generic filter over $V$.
Then in $V[G]$ every first countable $<\!\aleph_2$-cwH space is weakly
cwH. \endproclaim

\demo {Proof}
We have to show that in $V[G]$, if $\theta \ge \omega_2$ is a cardinal 
and $X$ is a first countable, $<\!\aleph_2$-cwH space, then $X$ is weakly 
$\theta$-cwH.
It is a consequence of a result of Watson \cite{W} that if $\theta$ is a
singular  strong limit cardinal, and $X$ is first countable and weakly
$<\!\theta$-cwH,  then $X$ is weakly $\theta$-cwH.  Therefore, we may assume
that  $\theta$ is regular. Since $\kappa = \omega_2$ in $V[G]$, we have
$\theta \ge \kappa$. The plan is to show that (2) of Theorem~2,  with
$\kappa = \omega_2$, fails in $V[G]$.

Let $j: V \to M$ be an elementary embedding such that
\roster
\item "(i)" $j(\alpha) = \alpha$, for all $\alpha < \kappa$.
\item "(ii)" $j(\kappa) > \theta$.
\item "(iii)" $^\theta M \cap V \subset M$.
\endroster

Now, ${\Bbb P} \subset V_\kappa$ and $j \restriction V_\kappa$ is the
identity map, therefore ${\Bbb P} \subset M$ and by (iii), ${\Bbb P}\in M$.
So, $M[G]$ makes sense.  
In $M$ (and also in $V$), $j({\Bbb P})$ is the Levy collapse of $j(\kappa)$
to $\omega_2$.  In the usual way $j({\Bbb P})$ can be factored, so   
$j({\Bbb P}) = j({\Bbb P})_\kappa \times j({\Bbb P})^\kappa$. But $j({\Bbb
P})_\kappa = {\Bbb P}$, so  $j({\Bbb P}) = {\Bbb P} \times j({\Bbb
P})^\kappa$. Moreover,
$$
M[G] \models  j({\Bbb P})^\kappa  \text{ is an $\omega_1$-closed
partial order.} 
$$

Let $I$ be a $j({\Bbb P})^\kappa$-generic filter over $V[G]$.
$I$ is also $j({\Bbb P})^\kappa$-generic over $M[G]$, and there is a filter
$K \subset j({\Bbb P})$ such that $M[G][I] = M[K]$ and $V[G][I] = V[K]$.
Using the fact that for
every $p \in {\Bbb P}$, $j(p) = p$, we can define, in $V[K]$, 
$$
j^* : V[G] \to M[K]
$$ 
that extends $j$ and show that it is an
elementary embedding (see \cite{KM} \S 25).

Let ${\Cal H} \in V[G]$ with 
${\Cal H}=\{H_{\alpha\beta} : \alpha < \beta < \theta\} 
\subset [\omega \times \omega]^{<\omega}$ c.d.w.\ and assume that 
the following condition from Theorem~2 holds in $V[G]$:
\roster
\item "(b)" For every $B \in [\theta]^{\theta}$ and every
$f:\theta \to \omega$ there are $\alpha <\beta$ in $B$
such that $(f(\alpha),f(\beta))\in H_{\alpha\beta}$.
\endroster
We will show that (a) from Theorem~2 (with $\kappa = \omega_2$) 
is false for ${\Cal H}$ in $V[G]$.

${\Cal H}$ has a name of size $\theta$ in $V$ that is a subset of $M$,
therefore by (iii) this name is in $M$.  Hence, ${\Cal H} \in M[G]$
and it is clear that condition (b) holds for ${\Cal H}$ in $M[G]$
as well. Notice that formally ${\Cal H}$ is a function given by 
${\Cal H} =\{( \langle \alpha , \beta \rangle, H_{\alpha\beta}) :  \alpha <
\beta < \theta\}$.  Let ${\Cal H'} = j^*({\Cal H}) \in M[K]$. Then
${\Cal H'}=\{\langle \langle \gamma , \delta \rangle, 
H'_{\gamma\delta}\rangle : \gamma < \delta < j(\theta)\}$, where each
$H'_{\gamma\delta}$ is a c.d.w. subset of 
$ [\omega \times \omega]^{<\omega}$.
Notice that for every $\alpha < \beta < \theta$,  $H_{\alpha\beta} =
j(H_{\alpha\beta}) = H'_{j(\alpha)j(\beta)}$. Let $A = j''\theta$, then 
$j''{\Cal H}  = {j^*}''{\Cal H} = \{( \langle
j(\alpha) , j(\beta) \rangle,  H_{\alpha\beta}) : \alpha < \beta <
\theta\} =  \{(\langle\gamma , \delta \rangle,  H'_{\gamma\delta}) : 
(\gamma < \delta)\land (\gamma, \delta \in A)\} \subset {\Cal H'}$.
The map $j \restriction {\Cal H}$ is a subset of $M$ of size
$\theta$, but it is not in $V$, therefore we can not conclude that it is
in $M$.  Fortunately,  
$^\theta M[G] \cap V[G] \subset M[G]$ 
(see \cite{J} p. 463), therefore 
$j \restriction {\Cal H} \in M[G] \subset M[K]$.

All that implies that the following two statements are equivalent in
$M[K]$:
\roster 
\item "(*)" for every $f : j(\theta) \to
\omega$ there are $\gamma < \delta$ in $A$ such that 
$(f(\gamma),f(\delta)) \in H'_{\gamma\delta}$.
\item "(**)" for every $f : \theta \to \omega$ there are $\alpha <
\beta < \theta $ such that  $(f(\alpha),f(\beta)) \in H_{\alpha\beta}$. 
\endroster 

As was metioned before, (b) holds for ${\Cal H}$ in $M[G]$ and
$M[G] \models ``j({\Bbb P})^\kappa$ is an $\omega_1$-closed
partial order".
Therefore, by the main lemma, $(**)$ (and hence $(*)$)
holds in  $M[K]$.

We conclude that in $M[K]$, there is $A \in
[j(\theta)]^{<j(\kappa)}$ such that for every $f: j(\theta) \to \omega$
there are $\gamma < \delta$ in $A$ with $(f(\gamma),f(\delta)) \in j^*({\Cal
H})_{\gamma \delta}\quad(= H'_{\gamma\delta})$.

Therefore, by elementarity of $j^*$, in $V[G]$ there is an
$A \in [\theta]^{\le \aleph_1}$ such that for every $f:\theta \to \omega$ 
there are $\alpha < \beta$ in $A$
with $(f(\alpha),f(\beta)) \in H_{\alpha \beta}$.
\qued

We can also prove Theorems~4 and~5 replacing the Levy collapse by the
Mitchell collapse (\cite{Mi} and \cite{DJW}), and thus obtain the
conclusions of these theorems in models where $CH$ fails.

\subhead \S 4 Integer-Valued Functions \endsubhead

As was mentioned before, Todorcevic proved that $\square(\theta)$ implies that
there is a  $\theta$-good subset of $\fsquare$~\cite{To2}.   In \cite{DW}, it is
proved that all one needs is a  non-extendible monotone $\theta$-family of
functions (which exists under $\square(\theta)$ but not in the model of
Theorem~5).

Let $g, h$ be functions with range a subset of $\omega$.  We say that 
{\it $g$ weakly bounds $h$\/} if there is $n \in \omega$ such that for
every $x \in dom(g) \cap dom(h)$ 
$$
g(x)+n > h(x)\qquad(\text{in short, }  g+n>h).
$$

We say that  $\Cal H = \{h_\beta : \beta < \theta \}$ is a $\theta$-family of
functions if the domain of $h_\beta$ is $\beta$.  The $\theta$-family $\Cal H$ is
{\it weakly bounded\/} if there is a $g$ that weakly bounds each element of $\Cal
H$.  We say that $\Cal H$ is an {\it initially weakly bounded
(i.w.b.) $\theta$-family\/}
 if for every $A \in [\theta]^{<\theta}$,  $\{h_\beta : \beta \in A \}$
is weakly bounded.  The family is {\it non-extendible\/} if $\{h_\beta : \beta <
\theta \}$ is not weakly bounded.

We can now repeat the construction of a $\theta$-good subset of
$\fsquare$ from a non-extendible i.w.b.
$\theta$-family $\{h_\beta : \beta < \theta \}$, using i.w.b. in place of
monotonicity.  To see this, define for every  $\alpha < \beta < \theta$ the
finite set $H_{\alpha\beta}=\{(n,m)\in \omega\times\omega :   n+m \le
h_\beta(\alpha)\}$ and check that ${\Cal H} = \{H_{\alpha\beta}: \alpha <
\beta < \theta\}$  satisfies~(3) of Theorem~1.

\proclaim {Question 3}
Does the existence of a $\theta$-good subset of $\fsquare$ imply the existence of a non-extendible i.w.b.
$\theta$-family of functions?
\endproclaim

Fleissner \cite{F} used the combinatorial principle $\etheta$ to
construct a locally countable, locally compact Moore space that is
$<\!\theta$-cwH but not $\le\!\theta$-cwH.  In particular, by \cite{GT} (see
$(1) \iff (2)$ of Theorem~1),  $\etheta$ implies
the existence of a $\theta$-good subset of $\fsquare$.

\proclaim {Theorem~6}
Let $\theta$ be an uncountable regular cardinal, and let  
$E \subset \{\alpha \in \theta : cof (\alpha) = \omega \}$ be a
non-reflecting stationary set.  
Then there is a $\theta$-family $\{h_\beta : \beta < \theta \}$
such that for every 
$A \subset \theta$, $\{h_\beta : \beta \in A \}$ is weakly
bounded if and only if $A \cap E$ is non-stationary. (In particular
$\{h_\beta : \beta < \theta \}$ is a non-extendible i.w.b. $\theta$-family.)
\endproclaim

\demo {Proof}
For every $\alpha \in E$, fix
$\{a_\alpha (n) : n \in \omega\}$ an increasing $\omega$-sequence unbounded
in $\alpha$.  For every $\alpha < \beta < \theta$ define
$$
h_\beta(\alpha) = \cases 
\min \{n \in \omega: a_\alpha(n) \ne a_\beta(n)\} 
&\text{if } (\alpha \in E) \land (\beta \in E), \\
0 & \text {otherwise}.
\endcases
$$

It is enough to consider $A \subset E$.  Let $A$ be stationary and assume, by
way of contradiction, that $g : \theta \to \omega$ is a weak bound for
$\{h_\beta : \beta \in A \}$.  There exists a stationary $B \subset A$ and
$n \in \omega$ such that for every $\beta \in B$, $g + n > h_\beta$.  
Let $f=g+n$, then for every $\beta \in B$, $f > h_\beta$.

Now, there exists a stationary $S \subset B$, and $m \in \omega$ such that
for every $\alpha, \beta \in S$, $f(\alpha) = m$ and 
$\langle a_\alpha (0), \dots, a_\alpha (m)\rangle = 
\langle a_\beta (0), \dots, a_\beta (m)\rangle$. 
Let $\alpha < \beta$ in $S$.  Let 
$k = \min \{n \in \omega: a_\alpha(n) \ne a_\beta(n)\}$.
Then $k > m$, therefore $f(\alpha) = m < k = h_\beta(\alpha)$, contradiction.

Next we prove, by induction on $\gamma < \theta$, that there exists 
$g_\gamma : \gamma \to \omega$ that weakly bounds $\{h_\beta : \beta<\gamma \}$.

The successor case is easy so let us assume that $\gamma$ is a limit.
let $C \subset \gamma$ be a club with $C \cap E \cap \gamma = \emptyset$.
For every $\beta \in E \cap \gamma$ let
$$
\align
&\gamma^+(\beta)=\min(C \setminus \beta) \\
&\gamma^-(\beta)=\max(C \cap \beta).
\endalign
$$

Notice that for every $\beta \in E \cap \gamma$ we have
$\gamma^-(\beta) < \beta < \gamma^+(\beta) < \gamma$.  Now define
$$
g_\gamma(\beta) = \cases 
g_{\gamma^+\!(\beta)}(\beta) &\text{if } \beta \in E, \\
1 & \text {otherwise}.
\endcases
$$

Let us show that $g_\gamma$ weakly bounds $\{h_\beta : \beta<\gamma \}$. If
$\beta\notin E$, then $h_\beta$ is identically zero, so it suffices to
consider $\beta \in E \cap \gamma$.  Let $n = k + m + 1$, where $k$ is the
least integer such that $a_\beta(k) > \gamma^-(\beta)$ and $m$ is satisfies
 $g_{\gamma^+\!(\beta)} + m > h_\beta$. Let us show that $g_\gamma + n >
h_\beta$. Let $\alpha \in E \cap \beta$.

\roster
\item "Case 1:"  $\gamma^-(\alpha)=\gamma^-(\beta)$.  In this case, 
$\gamma^+(\alpha)=\gamma^+(\beta)$, therefore
$g_\gamma(\alpha) = g_{\gamma^+\!(\alpha)}(\alpha) =
g_{\gamma^+\!(\beta)}(\alpha)$. But $g_{\gamma^+\!(\beta)}(\alpha) + m >
h_\beta(\alpha)$, therefore $g_\gamma(\alpha) + n > h_\beta(\alpha)$.

\item "Case 2:"  $\alpha < \gamma^-(\beta)$.  Then  $k+1 > h_\beta(\alpha)$
because $a_\beta(k) > \gamma^-(\beta)$, therefore $g_\gamma(\alpha) + n >
h_\beta(\alpha)$. \endroster

Finally, let $A \subset E$ be non-stationary.  Let us use  $\{g_\gamma :
\gamma < \theta\}$ to produce $g : \theta \to \omega$ that weakly bounds
$\{h_\beta : \beta \in A \}$.  Let $C \subset \theta$ be a club such that $C
\cap A = \emptyset$.  For every $\beta \in E$ let 
$$
\align
&\delta^+(\beta)=\min(C \setminus (\beta + 1)) \\
&\delta^-(\beta)=\sup(C \cap \beta), \\
&g(\beta) = \cases 
g_{\delta^+\!(\beta)}(\beta) & \text{if }\beta \in E, \\
1 & \text{otherwise.}
\endcases
\endalign
$$

Let $\beta \in A$.  As before, let $n = k + m + 1$, where 
$k$ is the least such
that $a_\beta(k) > \delta^-(\beta)$ (notice, $\beta \in A$ implies that 
$\delta^-(\beta) < \beta$),
and $m$ is such that 
$g_{\delta^{+}\!(\beta)} + m > h_\beta$.
Let us show that $g + n > h_\beta$.
Let $\alpha \in E \cap \beta$.

\roster
\item "Case 1:"  $\delta^-(\alpha)=\delta^-(\beta)$.  So, 
$\delta^+(\alpha)=\delta^+(\beta)$, therefore
$g(\alpha) = g_{\delta^+\!(\alpha)}(\alpha) =
g_{\delta^+\!(\beta)}(\alpha)$. But $g_{\delta^+\!(\beta)}(\alpha) + m >
h_\beta(\alpha)$, therefore $g(\alpha) + n > h_\beta(\alpha)$.
\item "Case 2:"  $\delta^-(\alpha)<\delta^-(\beta)$.  In this case 
$\delta^-(\beta) \ge \alpha$.  But $a_\beta(k) > \delta^-(\beta)$,
therefore $a_\beta(k) > \alpha$ and hence 
$a_\beta(k) \ne a_\alpha(k)$.  So, $k+1 > h_\beta(\alpha)$, therefore
$g(\alpha) + n > h_\beta(\alpha)$.
\endroster
\qued

We remark here that $\neg \eset$ is equiconsistent with a Mahlo
cardinal (\cite{D} and \cite{HS}), while  $\neg\square(\omega_2)$ is 
equiconsistent with a weakly compact cardinal (\cite{Ma}, \cite{To1} and
\cite{Be}),
(and therefore $\square(\omega_2)$ does not imply $\eset$).   And so
consistency-wise, a result from $\square(\omega_2)$ is better (more
difficult to deny) than one from $\eset$.  On the other hand,  $\eset$ does
not imply $\square(\omega_2)$ (PFA is consistent with $\eset$ but implies
$\neg \square(\omega_2)$ \cite{Be}).

\proclaim {Theorem~7}
Let $\kappa$ be a supercompact cardinal, let ${\Bbb P}$ be the Levy
collapse of $\kappa$ to $\omega_2$ with countable conditions and let $G$ be
a ${\Bbb P}$-generic filter over $V$.
Then in $V[G]$ the following holds:
For every regular $\theta \ge \omega_2$, and every
$\theta$-family $\{h_\beta : \beta < \theta \}$, 
if for all $A \in [\theta]^{\aleph_1}$,  $\{h_\beta : \beta \in A \}$
is weakly bounded, then for every stationary
set $S \subset \theta$, there exists a stationary $T \subset S$
such that $\{h_\beta : \beta \in T \}$ is weakly bounded.
\endproclaim

The proof of Theorem~7 is very much like th proof of Theorem~5.  The
following lemma is the analog of the Main Lemma, the
rest of the proof of the theorem is left to the reader.  

\proclaim {Lemma}
Assume that ${\Bbb P}$ is an $\omega_1$-closed partial order, and that 
$\theta \ge \omega_2$ is a regular cardinal.  Assume that 
${\Cal H} = \{h_\beta : \beta < \theta \}$ is a $\theta$-family, and 
$S \subset \theta$
is stationary such that for every stationary $T \subset S$ the family
$\{h_\beta : \beta \in T \}$ is not weakly bounded.
Assume that $G$ is a ${\Bbb P}$-generic filter over $V$. Then in
$V[G]$, $\{h_\beta : \beta \in S \}$ is not weakly bounded.
\endproclaim

\demo {Proof}
We proceed as in our Main Lemma and assume the lemma is false.  Let 
$p_0 \in {\Bbb P}$ and $\dot g$, a ${\Bbb P}$-name, be such that 
$p_0 \forces  ``\dot g: \theta \to \omega$ weakly bounds 
$\{h_\beta : \beta \in S \}"$.

Let $M$ be an elementary substructure of $H(\lambda)$ for some large enough
regular $\lambda$, with  $S$, $\theta$, $\Cal H$, $\Bbb P$, $p_0$, $\dot g$, and
$\forces_{\Bbb P}$
all members of $M$ and  $M \cap \theta = \beta \in S$.

Construct $p_0 \ge p_1 \ge \dots \ge p_{n+1} \dots$ such that
$p_{n+1} \in M$ and $\alpha_n < \beta$ with
$$
p_{n+1} \forces `` \dot g(\alpha_n) + n \le h_\beta(\alpha_n)".
$$
This process must stop at a finite stage since otherwise one could choose 
$q$ with $q \le p_n$ for every $n \in \omega$, and get that
for every $n \in \omega$ there is $\alpha_n < \beta$ such that
$q \forces `` \dot g(\alpha_n) + n \le h_\beta(\alpha_n)"$, 
which is impossible.

So, there is $\bar n \in \omega$ such that for every $q \le p_{\bar n}$ and
for every $\alpha < \beta$, if $q \in M$, then 
$q \nforces ``\dot g(\alpha) + \bar n \le  h_\beta(\alpha)"$.
As in the Main Lemma, it follows that 
for every $\alpha < \beta$,
$$ 
p_{\bar n} \forces ``\dot g(\alpha) + \bar n >  h_\beta(\alpha)".
$$ 
Let $T=\{ \beta \in S :  (\forall \alpha < \beta) \,\,
p_{\bar n} \forces ``\dot g(\alpha) + \bar n >  h_\beta(\alpha)"\}$. 
By elementarity of $M$, $T$ is stationary.  

Now, for every $\alpha \in \theta$, let $q_\alpha \le p_{\bar n}$ and 
$m_\alpha \in \omega$ be such that $q_\alpha \forces `` \dot 
g(\alpha) = m_\alpha"$.
 In $V$, define $f : \theta \to \omega$ by 
$f(\alpha) = m_\alpha + {\bar n}$.  Finally, continue as in the Main Lemma
to get a contradiction by showing that for every $\beta \in T$
and every $\alpha < \beta$, $   h_\beta(\alpha)<f(\alpha)$.
In particular, $f$ weakly bounds $\{h_\beta : \beta \in T \}$.
\qued

By Theorem~6, in the model of Theorem~7, $\etheta$ fails for every regular
$\theta \ge \omega_2$, (this is a result of Shelah).  As was the case for
Theorems~4 and~5, Theorem~7 can also be proved for the Mitchell collapse
instead of the Levy collapse.  We finish with a question.

\proclaim {Question~4}
Is it consistent (relative to a large cardinal) that every
i.w.b.  $\omega_2$-family extends? \endproclaim

\enddocument